# Why do launch trajectories end downwards

Max CERF [*]
Ariane Group, France

## Abstract

The problem of finding the optimal thrust profile of a launcher upper stage is analyzed. The engine is assumed to be continuously thrusting, following either a linear or a bilevel parametric profile, until reaching the targeted coplanar orbit. The minimum-fuel problem is analytically solved using the maximum principle. A closed-loop solution for the thrust direction is derived and the final point is found to be necessarily at an apside, reached downwards in the case of a perigee injection. The optimal control problem reduces to a nonlinear problem with only the thrust profile parameters as unknowns. An application case targeting a geostationary transfer orbit illustrates the solution method.

## 1. Introduction

The thrust levels of the engines are key parameters of a launcher design. They drive the loads met during the flight, the fuel required to reach the targeted orbit, and finally the launcher gross mass and cost. The problem of finding the best thrust levels is thus intensively investigated since Goddard pioneering work [1-4].

For a given launcher, the thrust levels are prescribed and they can generally not be changed. The ascent trajectory is controlled through the thrust orientation. Finding the minimum-fuel trajectory reaching the orbit with path constraints is a classical optimal control problem [5-13]. Many efforts have been devoted to theoretically analyze this problem. Exact or approximate analytical solutions may be derived depending on the dynamics formulation : impulsive maneuvers [14-16], high thrust [11,13], low thrust [6,7,13,17,18], endo-atmospheric flight [19-21], exo-atmospheric flight [22,23], constant gravity [9,23], linear gravity [22], ...
Such analytical solutions, when available, are of limited scope. In order to solve real-world applications a variety of numerical methods have been developed that divide into direct and indirect methods.

---

[*] Ariane Group, 78130 Les Mureaux, France
max.cerf@airbusafran-launchers.com

Direct methods [24-32] transform the optimal control problem in a nonlinear programming problem. Many discretizing approaches can be envisioned. The resulting large-scale optimization problem is then solved by a nonlinear (NLP) software (IPOPT, BOCOP, GESOP, SNOPT, WORHP …). Direct methods handle easily any problem formulation (dynamics, constraints) with reduced programming effort. Due to the problem large size they may nevertheless be computationally expensive and possibly inaccurate.

Indirect methods [33-36] transform the optimal control problem into a boundary value problem. The optimal control is explicitly or implicitly determined from the Pontryaguin Maximum Principle (PMP). The problem unknowns are the initial costates that must be found in order to satisfy transversality conditions at the final time. The small problem size makes the indirect approach attractive, but numerical issues arise due to the high sensitivity to the initial guess, and possibly to control discontinuities along the trajectory [4,37-39].

Hybrid approaches strive to combine direct and indirect methods in order to benefit from their respective qualities [40,41]. A direct method is first used to build a good initial guess, and an indirect method is applied to yield an accurate convergence. Amongst the numerous hybrid techniques applied to trajectory optimization, we can mention impulsive solution guess [42-45], multiple shooting [44], homotopy [46,47], response surface [48], averaging [49-51], particle swarm [52], genetic algorithms [53], dynamic programming based on the Hamilton-Jacobi-Bellman [54], differential dynamic programming [55-57], …

The above methods deal mainly with pure trajectory problems considering a given vehicle configuration. Preliminary design studies aim at finding the optimal launcher configuration for a reference mission. This results in a multidisciplinary problem addressing together configuration parameters and trajectory optimization. Due to their overall complexity such multidisciplinary problems cannot be addressed in a one-shot optimization process. A general procedure consists in iterating between configuration changes and trajectory optimization, until convergence toward an acceptable design. The iterative process can be built in many ways depending on the vehicle optimized sub-systems [53,58-62]. Among other configuration parameters the engine thrust levels play a key role in the launcher design and performance. More specifically the thrust levels of the lower stages are driven by loads considerations (acceleration, dynamic pressure, thermal flux, …), while the upper stage thrust level is driven by fuel minimization. This last parameter is thus the most strongly coupled to the trajectory optimization problem. Solving efficiently this last stage design problem is thus of utmost interest.

This paper addresses the minimum-fuel single-boost planar trajectory of a launcher upper stage. The initial conditions are sub-orbital resulting from the previous stage flight. The targeted coplanar orbit has a defined shape and a free orientation. This planar problem is quite close to most practical applications. The thrust level profile is part of the optimization problem. Two parametric models are considered, either linear or bilevel. Instead of complicating the optimal control problem, these additional parameters allow deriving a closed-loop control solution for the thrust direction. The optimal injection point is also found to be necessarily at an apside of the targeted orbit. These theoretical properties reduce the optimal control problem to a small size nonlinear problem which can be easily solved, yielding the upper stage optimal thrust profile and trajectory for a given mission.

The text is divided into a section §2 formulating the optimal control problem and deriving the solution properties, and a section §3 presenting numerical applications.

## 2. Theoretical Analysis

The Optimal Control Problem (OCP) is formulated and the Pontryaguin Maximum Principle (PMP) is applied. A closed-loop control is derived for the thrust direction, and the injection point is found to be at an apside.

### 2.1 Optimal Control Problem

The problem consists in finding the planar minimal-fuel trajectory going from given injection conditions to a given orbit. The dynamics equations are written in an Earth centered inertial frame.

$$\begin{cases} \dot{\vec{r}}(t) = \vec{v}(t) \\ \dot{\vec{v}}(t) = \vec{g}[\vec{r}(t)] + \dfrac{T(t)}{m(t)} \vec{u}(t) \\ \dot{m}(t) = -\dfrac{T(t)}{v_e} \end{cases} \tag{1}$$

where  $\vec{r}(t)$, $\vec{v}(t)$, $m(t)$ denote respectively the vehicle position, velocity and mass at the date t

$\vec{g}(\vec{r}) = -\dfrac{\mu}{r^3} \vec{r}$ is the Earth acceleration gravity at the position $\vec{r}$, with the gravitational constant $\mu$

$T(t)$ is the thrust level at the date t

$\vec{u}(t)$ is a unit vector defining the thrust direction at the date t

$v_e$ is the burned propellant exhaust velocity.

The trajectory is planar and all vectors are of dimension two.

The upper stage is released at the date $t_0$ on a sub-orbital trajectory corresponding to the fall-out of the last but one stage. The initial conditions at $t_0$ are denoted $\vec{r}_0$, $\vec{v}_0$ and $m_0$. The engine is ignited at the initial date $t_0$ and it will then thrust constantly until reaching the targeted orbit at a date $t_f$.

The targeted coplanar orbit is defined either by apogee and perigee altitudes, or semi major axis and eccentricity, or energy and angular momentum modulus per mass unit. All these definitions are equivalent for a planar transfer. The orbit orientation is free. The energy and angular momentum are given by

$$\begin{cases} w = \dfrac{v^2}{2} - \dfrac{\mu}{r} = -\dfrac{\mu}{2a} \\ \vec{h} = \vec{r} \times \vec{v} \quad, \quad h = \sqrt{\mu a (1-e^2)} \end{cases} \tag{2}$$

where  a, e denote respectively the semi-major axis and the eccentricity

w is the energy per unit mass

$\vec{h}(t)$ is the angular momentum per unit mass, orthogonal to the transfer plane

The targeted values of the energy and angular momentum modulus are denoted respectively $w_f$ and $h_f$.

Two parametric models are considered for the thrust level evolution, respectively a linear and a bilevel model.

- Linear model

$$T(t) = T_1 + (t-t_0)T_2 \quad \text{for } t_0 \leq t \leq t_f \tag{3}$$

- Bilevel model

$$T(t) = \begin{cases} T_1 & \text{if } t_0 \leq t < t_1 \\ T_2 & \text{if } t_1 \leq t < t_f \end{cases} \tag{4}$$

The thrust level T depends therefore on the two parameters $T_1$ and $T_2$, and also on the switching date $t_1$ for the bilevel model. It is assumed that the thrust remains strictly positive on $[t_0, t_f]$.

The optimal control problem (OCP) consists in finding the values of $T_1$, $T_2$, $t_1$, $t_f$ together with the thrust direction $\vec{u}(t)$, $t_0 \leq t \leq t_f$ in order to maximize the final mass injected on the targeted orbit.

Optimal Control Problem (OCP)

$$\min_{\substack{T_1, T_2, t_1, t_f \\ \vec{u}(t)}} J = -m(t_f) \quad \text{s.t.} \quad \begin{cases} \dot{\vec{r}} = \vec{v} \\ \dot{\vec{v}} = \vec{g} + \dfrac{T}{m}\vec{u} \\ \dot{m} = -\dfrac{T}{v_e} \end{cases} \text{with} \quad \begin{cases} \vec{r}(t_0) = \vec{r}_0 \\ \vec{v}(t_0) = \vec{v}_0 \quad \text{fixed initial state} \\ m(t_0) = m_0 \end{cases} \quad \begin{cases} w(t_f) = w_f \\ h(t_f) = h_f \end{cases} \text{final constraints} \tag{5}$$

The Hamiltonian for this optimal control problem is

$$\begin{aligned} H &= \vec{p}_r^t \dot{\vec{r}} + \vec{p}_v^t \dot{\vec{v}} + p_m \dot{m} = \vec{p}_r^t \vec{v} + \vec{p}_v^t \left(\vec{g} + \dfrac{T}{m}\vec{u}\right) + p_m \left(-\dfrac{T}{v_e}\right) \\ &= \vec{p}_r^t \vec{v} + \vec{p}_v^t \vec{g} + T\left(\dfrac{\vec{p}_v^t \vec{u}}{m} - \dfrac{p_m}{v_e}\right) \end{aligned} \tag{6}$$

where $\vec{p}_r(t)$, $\vec{p}_v(t)$, $p_m(t)$ are the costate vectors respectively associated to the position, the velocity, the mass. These vectors do not vanish identically on any interval of $[t_0, t_f]$ and they are defined up to a non-positive scalar multiplier $p_0$ which can be chosen freely for a normalization purpose.

The Pontryaguin Maximum Principle (PMP) [33-35] provides first order necessary optimality conditions that must be satisfied along an optimal trajectory.

- The Hamiltonian maximization condition yields the thrust direction $\vec{u}(t)$.

$$\max_{\vec{u}} H \quad \Rightarrow \quad \vec{u} = \dfrac{\vec{p}_v}{p_v} \quad \text{with} \quad p_v = \|\vec{p}_v\| \tag{7}$$

Replacing $\vec{u}(t)$ in Eq. (6), the Hamiltonian takes the form

$$H = H_0 + T\Phi \quad \text{with} \quad H_0 \underset{\text{def}}{=} \vec{p}_r^t \vec{v} + \vec{p}_v^t \vec{g} \quad \text{and} \quad \Phi \underset{\text{def}}{=} \dfrac{p_v}{m} - \dfrac{p_m}{v_e} \tag{8}$$

- The costate differential equations are

$$\begin{cases} \dot{\vec{p}}_r = -\dfrac{\partial H}{\partial \vec{r}} = -\dfrac{\partial \vec{g}}{\partial \vec{r}} \vec{p}_v \\ \dot{\vec{p}}_v = -\dfrac{\partial H}{\partial \vec{v}} = -\vec{p}_r \\ \dot{p}_m = -\dfrac{\partial H}{\partial m} = \dfrac{T}{m^2} \vec{p}_v^T \vec{u} = \dfrac{T}{m^2} p_v \quad \text{using} \quad \vec{u} = \dfrac{\vec{p}_v}{p_v} \quad \text{from Eq. (7)} \end{cases} \qquad (9)$$

- The final costates satisfy the transversality conditions coming from the final constraints on w and h with respective multipliers $\nu_1$, $\nu_2$ and from the final cost $-m(t_f)$ with the non-positive multiplier $p_0$.

$$\begin{cases} \vec{p}_r(t_f) = \nu_1 \dfrac{dw}{d\vec{r}}(t_f) + \nu_2 \dfrac{dh}{d\vec{r}}(t_f) \\ \vec{p}_v(t_f) = \nu_1 \dfrac{dw}{d\vec{v}}(t_f) + \nu_2 \dfrac{dh}{d\vec{v}}(t_f) \\ p_m(t_f) = -p_0 \end{cases} \qquad (10)$$

- The Hamiltonian is constant since the problem is autonomous. Its final value comes from the transversality condition associated to the free final date $t_f$.

$$H(t) = H(t_f) = 0 \, , \quad \forall t \in [t_0, t_f] \qquad (11)$$

- The parameters $T_1$ and $T_2$ must satisfy the following first order conditions [34].

$$\begin{aligned} \int_{t_0}^{t_f} \dfrac{\partial H}{\partial T_1} dt = 0 &\Rightarrow \int_{t_0}^{t_f} \dfrac{\partial T}{\partial T_1} \Phi \, dt = 0 \\ \int_{t_0}^{t_f} \dfrac{\partial H}{\partial T_2} dt = 0 &\Rightarrow \int_{t_0}^{t_f} \dfrac{\partial T}{\partial T_2} \Phi \, dt = 0 \end{aligned} \qquad \text{using} \quad H = H_0 + T\Phi \quad \text{from Eq. (8)} \qquad (12)$$

- The date $t_1$ define a dynamics change at an interior point. If a constraint applies at such a point, it generates discontinuities on the Hamiltonian and the costates [6]. In the present case, there are no such constraints. The Hamiltonian and the costates are continuous at $t_1$.

These first order optimality conditions can now be used, on the one hand to derive a closed-loop solution for the thrust direction (§2.2), on the other hand to find the optimal injection point (§2.3).

## 2.2 Thrust Direction

The integral conditions Eq. (12) depend on the function $\Phi$ and on the thrust level derivatives.

The function $\Phi$ is first rewritten as

$$\begin{aligned} \Phi = \dfrac{p_v}{m} - \dfrac{p_m}{v_e} &= \dfrac{m\dot{p}_m}{T} - \dfrac{p_m}{v_e} \qquad &\text{using} \quad \dot{p}_m = \dfrac{T}{m^2} p_v \quad \text{from Eq. (9)} \\ &= \dfrac{m\dot{p}_m + \dot{m} p_m}{T} \qquad &\text{using} \quad \dot{m} = -\dfrac{T}{v_e} \quad \text{from Eq. (1)} \\ &= \dfrac{\dot{\Psi}}{T} \qquad &\text{with} \quad \Psi \underset{\text{def}}{=} m p_m \end{aligned} \qquad (13)$$

with the assumption that the thrust does not vanish on $[t_0, t_f]$.

For the linear model Eq. (3), the thrust level derivatives are

$$\frac{\partial T}{\partial T_1} = 1 \quad \text{and} \quad \frac{\partial T}{\partial T_2} = t - t_0 \tag{14}$$

The conditions Eq. (12) take the form

$$\int_{t_0}^{t_f} \frac{\dot{\Psi}}{T} dt = 0 \quad \text{and} \quad \int_{t_0}^{t_f} \frac{\dot{\Psi}}{T}(t - t_0) dt = 0 \tag{15}$$

Integrating by parts the first condition

$$\int_{t_0}^{t_f} \frac{\dot{\Psi}}{T} dt = \left[\frac{\Psi}{T}\right]_{t_0}^{t_f} + \int_{t_0}^{t_f} \Psi \frac{\dot{T}}{T^2} dt = \left[\frac{\Psi}{T}\right]_{t_0}^{t_f} + \int_{t_0}^{t_f} \Psi \frac{T_2}{T^2} dt = 0 \tag{16}$$

and the second condition

$$\int_{t_0}^{t_f} \frac{\dot{\Psi}}{T}(t-t_0) dt = \left[\frac{\Psi}{T}(t-t_0)\right]_{t_0}^{t_f} - \int_{t_0}^{t_f} \Psi \frac{T - (t-t_0)\dot{T}}{T^2} dt = \left[\frac{\Psi}{T}(t-t_0)\right]_{t_0}^{t_f} - \int_{t_0}^{t_f} \Psi \frac{T_1}{T^2} dt = 0 \tag{17}$$

where Eq. (3) has been used to reduce the numerator in the integral.

The integral term is eliminated by adding Eq. (16) multiplied by $T_1$ to Eq. (17) multiplied by $T_2$.

$$T_1 \left[\frac{\Psi}{T}\right]_{t_0}^{t_f} + T_2 \left[\frac{\Psi}{T}(t-t_0)\right]_{t_0}^{t_f} = 0 \tag{18}$$

Grouping the terms and using again Eq. (3), this equation reduces to

$$[\Psi]_{t_0}^{t_f} = \Psi(t_f) - \Psi(t_0) = 0 \tag{19}$$

Coming back to the definition of $\Psi$ from Eq. (13), with $p_m(t_f)$ given by the transversality condition Eq. (10)

$$p_m(t_0)m(t_0) = -p_0 m(t_f) \tag{20}$$

For the bilevel model Eq. (4), the thrust level derivatives are

$$\frac{\partial T}{\partial T_1} = \begin{cases} 1 & \text{if } t < t_1 \\ 0 & \text{if } t > t_1 \end{cases} \quad \text{and} \quad \frac{\partial T}{\partial T_2} = \begin{cases} 0 & \text{if } t < t_1 \\ 1 & \text{if } t > t_1 \end{cases} \tag{21}$$

The conditions Eq. (12) take the form

$$\int_{t_0}^{t_1} \frac{\dot{\Psi}}{T} dt = 0 \quad \text{and} \quad \int_{t_1}^{t_f} \frac{\dot{\Psi}}{T} dt = 0 \tag{22}$$

T being piecewise constant, respectively on the intervals $[t_0, t_1]$ and $[t_1, t_f]$, we obtain

$$\begin{cases} [\Psi]_{t_0}^{t_1} = \Psi(t_1^-) - \Psi(t_0) = 0 \\ [\Psi]_{t_1}^{t_f} = \Psi(t_f) - \Psi(t_1^+) = 0 \end{cases} \tag{23}$$

No interior constraint applies at the switching date $t_1$, so that the function $\Psi = m p_m$ is continuous.

$$p_m(t_0)m(t_0) = p_m(t_1)m(t_1) = p_m(t_f)m(t_f) = -p_0 m(t_f) \tag{24}$$

The same condition Eq. (20) holds thus for the linear and the bilevel thrust model. The initial date $t_0$ being arbitrary this holds at any date t along the optimal trajectory.

$$p_m(t)m(t) = -p_0 m(t_f) \quad , \quad \forall t \in [t_0, t_f] \tag{25}$$

Derivating $p_m$ given by Eq. (25) and equating to $\dot{p}_m$ given by Eq. (9), we obtain an expression for $p_v$.

$$p_v(t) = -\frac{p_0 m(t_f)}{v_e} \quad , \quad \forall t \in [t_0, t_f] \tag{26}$$

The modulus $p_v$ of the velocity costate is constant along the optimal trajectory. Closed-loop solutions for the thrust direction and the costate vectors can be derived from this property. The detailled calculations are presented in [63] and they are not reproduced here. The main results are recalled here below.

The Figure 1 depicts the vectors and angles used. $\varphi$ is the longitude measured from the x axis in the Earth-centered Galilean frame (O,x,y). $\theta$ and $\gamma$ are respectively the local pitch angle and the flight path angle measured positively from the local horizontal. $\omega$ is the angular rate of $\vec{u}$ in the Galilean frame (O,x,y).

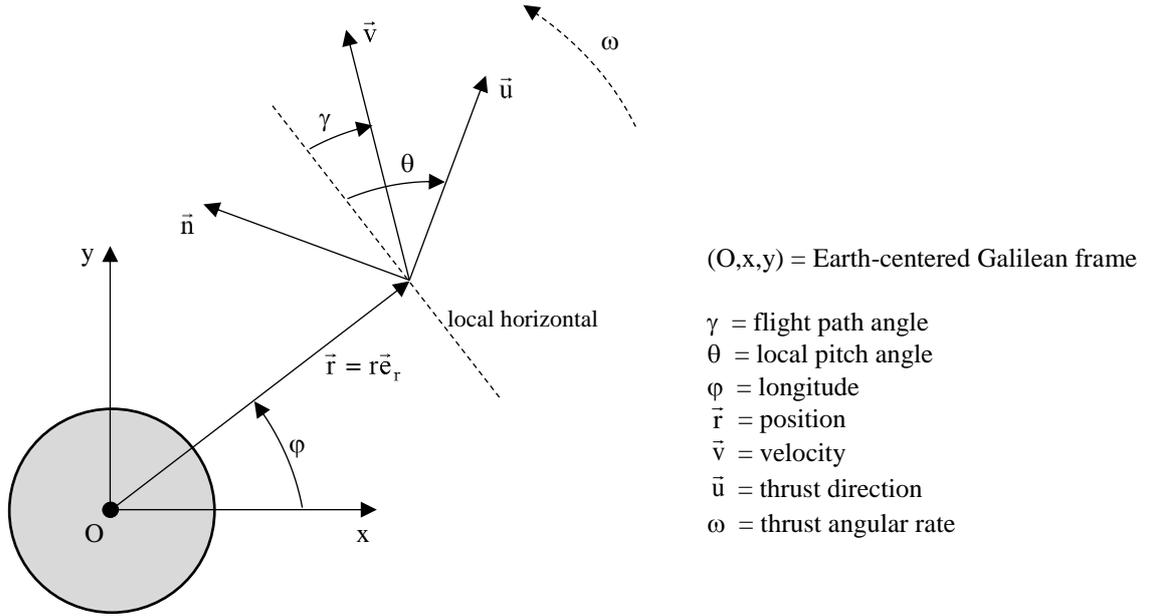

(O,x,y) = Earth-centered Galilean frame

$\gamma$ = flight path angle
$\theta$ = local pitch angle
$\varphi$ = longitude
$\vec{r}$ = position
$\vec{v}$ = velocity
$\vec{u}$ = thrust direction
$\omega$ = thrust angular rate

Figure 1 : Vectors and angles

With these notations the following conditions hold along the optimal trajectory.

- The thrust direction $\theta$ depends on the current kinematic conditions (r, v, $\gamma$) through

$$\sin(\gamma - \theta) = \frac{v_c}{v} \frac{\sin\theta}{\sqrt{1 - 3\sin^2\theta}} \quad \text{with} \quad v_c = \sqrt{\frac{\mu}{r}} \tag{27}$$

- The angular rate ω and its derivative are linked to θ, r, v, γ through

$$\begin{cases} \omega = \sqrt{\dfrac{\mu}{r^3}(1 - 3\sin^2\theta)} \\ \omega r v \sin(\theta - \gamma) = v_c^2 \sin\theta \\ \dot{\omega} = -\dfrac{3\mu}{r^3}\sin\theta\cos\theta \end{cases} \quad (28)$$

- The costates $\vec{p}_r(t)$, $\vec{p}_v(t)$ are given up to a constant multiplier by

$$\vec{p}_r = \omega \begin{pmatrix} -\cos(\theta-\varphi) \\ \sin(\theta-\varphi) \end{pmatrix} \quad \text{and} \quad \vec{p}_v = \begin{pmatrix} \sin(\theta-\varphi) \\ \cos(\theta-\varphi) \end{pmatrix} \quad (29)$$

Eq. (27) defines implicitly the thrust direction θ as a function of the kinematic conditions (see [63] for a detailled study of the possible solutions of this equation), so that $\vec{u}(t)$ is no longer an unknown. The OCP Eq. (5) reduces to a nonlinear programming problem (NLP) with $T_1$, $T_2$, $t_1$, $t_f$ as unknowns.

## 2.3 Injection Conditions

The transversality conditions Eq. (10) link the final costates to the constraint derivatives. These derivatives can be explicitly calculated for the energy and the angular momentum.

- The energy derivatives are directly assessed using the fact that the gravity derives from the potential.

$$w = \dfrac{v^2}{2} - \dfrac{\mu}{r} \quad \Rightarrow \quad \begin{cases} \dfrac{\partial w}{\partial \vec{r}} = -\vec{g} \\ \dfrac{\partial w}{\partial \vec{v}} = \vec{v} \end{cases} \quad \text{with} \quad \vec{g} = -\dfrac{\mu}{r^3}\vec{r} \quad (30)$$

- The angular momentum modulus is squared and halved, so that the derivatives can be easily assessed.

$$k \underset{\text{def}}{=} \dfrac{1}{2}h^2 = \dfrac{1}{2}(\vec{r}\times\vec{v})^t(\vec{r}\times\vec{v}) \quad \Rightarrow \quad \begin{cases} \dfrac{\partial k}{\partial \vec{r}} = -\vec{h}\times\vec{v} \\ \dfrac{\partial k}{\partial \vec{v}} = \vec{h}\times\vec{r} \end{cases} \quad (31)$$

The constraint on k is equivalent to the constraint on h. Replacing these derivatives in the transversality equations Eq. (10) yield two conditions at the final date $t_f$.

$$\begin{cases} \vec{p}_r = -\nu_1\vec{g} - \nu_2\vec{h}\times\vec{v} \\ \vec{p}_v = \nu_1\vec{v} + \nu_2\vec{h}\times\vec{r} \end{cases} \quad (32)$$

Taking the dot products respectively with $\vec{r}$ and $\vec{v}$, and using the equalities $(\vec{h}\times\vec{v})^t\vec{r} = \vec{h}^t(\vec{v}\times\vec{r}) = -\vec{h}^t\vec{h} = -h^2$ and $(\vec{h}\times\vec{r})^t\vec{v} = \vec{h}^t(\vec{r}\times\vec{v}) = \vec{h}^t\vec{h} = h^2$, we obtain three scalar equations.

$$\begin{cases} \vec{p}_r^t\vec{r} = \nu_1 v_c^2 - \nu_2 h^2 \\ \vec{p}_v^t\vec{v} = \nu_1 v^2 + \nu_2 h^2 \\ \vec{p}_v^t\vec{r} = \nu_1 \vec{r}^t\vec{v} \end{cases} \quad \text{with} \quad v_c = \sqrt{\dfrac{\mu}{r}} \quad (33)$$

The two first equations are summed to eliminate $\nu_2$.

$$\begin{cases} \vec{p}_r^t \vec{r} + \vec{p}_v^t \vec{v} = \nu_1(v^2 + v_c^2) \\ \vec{p}_v^t \vec{r} = \nu_1 \vec{r}^t \vec{v} \end{cases} \qquad (34)$$

The multiplier $\nu_1$ is eliminated by linear combination yielding a relation between the position, the velocity and their respective costates at the final date.

$$(\vec{p}_r^t \vec{r} + \vec{p}_v^t \vec{v}).(\vec{r}^t \vec{v}) = (\vec{p}_v^t \vec{r}).(v^2 + v_c^2) \qquad (35)$$

The components of the position and velocity vectors in the frame (O,x,y) are (see Figure 1)

$$\vec{r} = r \begin{pmatrix} \cos\varphi \\ \sin\varphi \end{pmatrix} \quad \text{and} \quad \vec{v} = v \begin{pmatrix} -\sin(\varphi-\gamma) \\ \cos(\varphi-\gamma) \end{pmatrix} \qquad (36)$$

while the costate components are given by Eq. (29), up to a multiplicative factor.
Replacing all vectors by their components in Eq. (35), we obtain a condition on the final conditions.

$$v^2 \cos\gamma \sin(\gamma-\theta) = v_c^2 \sin\theta + \omega r v \sin\gamma \cos\theta \qquad (37)$$

On the other hand, Eq. (28) provides another condition that is satisfied along the optimal trajectory.

$$\omega r v \sin(\theta-\gamma) = v_c^2 \sin\theta \qquad (38)$$

Using Eq. (38) the second member of Eq. (37) simplifies. Assuming that $v \neq 0$ and $\cos\gamma \neq 0$ (such conditions do not occur in practice), Eq. (37) becomes

$$v \sin(\gamma-\theta) = \omega r \sin\theta \qquad (39)$$

Multiplying both members by $\omega r$, and using again Eq. (38) we obtain the following condition at the final date.

$$(v_c^2 + \omega^2 r^2)\sin\theta = 0 \qquad (40)$$

This necessary condition can only be checked if the final thrust direction is horizontal. From Eq. (27) the final flight path angle is also necessarily null.

$$\sin\theta(t_f) = 0 \quad \Rightarrow \quad \gamma(t_f) = 0 \qquad (41)$$

<u>For an elliptical targeted orbit the optimal injection occurs either at the perigee or at the apogee.</u>

An additional information on the injection conditions comes from the derivative of $\theta$.
The thrust angular rate $\omega$ in the inertial frame is linked to the derivatives of $\varphi$ and $\theta$ (see Figure 1).

$$\omega = \dot\varphi - \dot\theta \quad \text{with} \quad \dot\varphi = \frac{v\cos\gamma}{r} \quad \text{and} \quad \omega = \sqrt{\frac{\mu}{r^3}(1-3\sin^2\theta)} \text{ from Eq.(28)} \qquad (42)$$

For $\theta = \gamma = 0$, this yields

$$\dot\theta = \frac{v - v_c}{r} \quad \text{with} \quad v_c = \sqrt{\frac{\mu}{r}} \qquad (43)$$

In the case of an injection at the perigee of an elliptical orbit, the velocity v is greater that the circular velocity $v_c$ and the derivative $\dot\theta$ is positive. The angles $\theta$ and $\gamma$ having the same sign (see Eq. (27) and Ref. [63] for the detailed study), the trajectory ends therefore with a negative flight path angle and the injection occurs downwards. Conversely an injection at the apogee occurs upwards.

In the case of a circular orbit, the second derivative must be assessed to find the angle sign.

$$\ddot{\theta} = \ddot{\varphi} - \dot{\omega} \quad \text{with} \quad \dot{\omega} = -\frac{3\mu}{r^3}\sin\theta\cos\theta \quad \text{from Eq.(28)} \tag{44}$$

For $\theta = \gamma = 0$, most terms vanish and we obtain

$$\ddot{\theta} = \frac{r\dot{v} - \dot{r}v}{r^2} = \frac{T}{mr} \quad \text{using} \quad \dot{r} = v\sin\gamma \quad \text{and} \quad \dot{v} = \frac{T}{m}\cos(\theta - \gamma) - g\sin\gamma \tag{45}$$

The thrust being strictly positive along the optimal trajectory, the second derivative is positive.
The trajectory ends with a negative flight path angle and the injection on the circular orbit occurs downwards.

In practice, most targeted orbits are either circular or with a high apogee. In such cases the optimal injection occurs at the perigee and the orbit is reached via a downward leg.

## 2.4 Performance Guess

The closed-loop control given by Eq. (27) can also be used to estimate the velocity losses along the optimal trajectory. The differential equation followed by the velocity modulus is

$$\dot{v} = \frac{T}{m}\cos(\theta - \gamma) - g\sin\gamma \tag{46}$$

The velocity losses come from the gravity and from the thrust unaligned with the velocity. These contributions called respectively gravity losses and angle of attack losses are denoted $\Delta V_G$ and $\Delta V_T$.

$$\begin{aligned} \Delta V_G &= \int_{t_0}^{t_f} g\sin\gamma\, dt \\ \Delta V_T &= \int_{t_0}^{t_f} \frac{T}{m}[1 - \cos(\theta - \gamma)]\, dt \end{aligned} \tag{47}$$

Along the optimal trajectory, the thrust angle $\theta$ satisfies Eq. (27)

$$\sin(\gamma - \theta) = \frac{v_c}{v}\frac{\sin\theta}{\sqrt{1 - 3\sin^2\theta}} \quad \text{with} \quad v_c = \sqrt{\frac{\mu}{r}} \tag{48}$$

For small values of $\gamma$ and $\theta$ (this occurs necessarily when approaching the horizontal injection), and considering that the velocities $v$ and $v_c$ have the same order of magnitude, this equation yields approximatively

$$\sin(\gamma - \theta) \approx \sin\theta \quad \Rightarrow \quad 2\theta \approx \gamma \tag{49}$$

On the other hand, the thrust angle is linked to the angular velocity derivative by Eq. (28).

$$\dot{\omega} = -\frac{3\mu}{r^3}\sin\theta\cos\theta \tag{50}$$

Combining Eqs. (49) and (50) we get

$$\dot{\omega} = -\frac{3\mu}{2r^3}\sin 2\theta \approx -\frac{3}{2r}g\sin\gamma \quad \text{with} \quad g = \frac{\mu}{r^2} \tag{51}$$

Replacing in Eq. (47), we get for the gravity losses.

$$\Delta V_G = \int_{t_0}^{t_f} g \sin \gamma \, dt \approx -\frac{2}{3} \int_{t_0}^{t_f} r \dot{\omega} \, dt \tag{52}$$

The radius vector variation is small along the upper stage propelled flight (typically between 100 and 300 km), so that it can be considered as constant, and the integral can be explicitly calculated.

$$\Delta V_G \approx \frac{2}{3} r [\omega_0 - \omega_f] \tag{53}$$

The angular rate is known at the trajectory endpoints owing to Eq. (28).

$$\omega = \sqrt{\frac{\mu}{r^3}(1 - 3\sin^2 \theta)} \tag{54}$$

At the initial date $t_0$, $r_0$ is given and $\theta_0$ is assessed either exactly from Eq. (27) or approximatively from Eq. (49). At the final date $t_f$, $r_f$ is known equal to $r_p$ (perigee injection) and $\theta_f$ is null.

Replacing in Eq. (54), assuming small flight path angles and a constant radius value equal to $r_0$, we obtain an estimate of the gravity losses along the optimal trajectory.

$$\Delta V_G \approx \frac{1}{4} v_c \gamma_0^2 \quad \text{with} \quad v_c = \sqrt{\frac{\mu}{r_0}} \tag{55}$$

This term $\Delta V_G$ is the main contributor to the velocity losses. From Eq. (49) the angles $\gamma$ and $\theta$ have indeed close values along the optimal trajectory, and the angle of attack losses $\Delta V_T$ remain small.

The total velocity impulse $\Delta V$ is the sum of the actual velocity increment $v_f - v_0$ and of the velocity losses. The final velocity is known since the injection occurs at an apside of the targeted orbit, generally at the perigee.

The final mass can then be estimated from the Tsiolkovsky formula (also called rocket equation).

$$\Delta V = v_e \ln \frac{m_0}{m_f} \quad \Rightarrow \quad m_f \approx m_0 e^{-\frac{\Delta V_G + v_p - v_0}{v_e}} \tag{56}$$

with $v_p$ the velocity at the perigee of the targeted orbit, and $\Delta V_G$ assessed with Eq. (55).

## 3. Application

The previous theoretical results help reducing the optimal control problem to a small size nonlinear programming problem. The optimal thrust direction $\theta$ depends only on the current kinematic conditions (r, v, γ) through Eq. (27). Solving this equation "on-line" yields a closed-loop control along the trajectory. The motion differential equations Eq. (1) are propagated using this closed-loop control until reaching a null flight path angle. For a perigee injection the final flight path angle must first decrease until a negative minimum, then increase until zero. This final condition determines the final time $t_f$. The remaining problem unknowns are the thrust level parameters $T_1$ and $T_2$, and optionally the switching time $t_1$ for the bilevel model. The constraints are the apogee and the perigee altitudes of the targeted orbit. The optimal control problem Eq. (1) reduces thus to a nonlinear problem with two or three unknowns. This problem is easily solved by any NLP software, starting from a rough initial guess.

The following example illustrates the solution method for a launcher upper stage targeting a coplanar geostationary transfer orbit (GTO). The targeted apogee and perigee have respective altitudes of 36000 km and 300 km. The engine specific impulse is 300 s corresponding to an exhaust velocity of 2942 m/s. The initial conditions at the stage ignition are a mass of 10000 kg, an altitude of 150 km, a velocity of 5000 m/s and a flight path angle of 30 deg, giving an apogee at 661.7 km and a perigee at −5209.4 km. The Earth equatorial radius is $R_E = 6378137$ m and the gravitational constant is $\mu = 3.986005 \cdot 10^{14}$ m$^3$/s$^2$.

A preliminary performance guess is first derived using Eqs. (55,56).

$$\begin{cases} m_0 = 10000 \text{ kg} \\ h_0 = 150 \text{ km} \\ v_0 = 5000 \text{ m/s} \\ \gamma_0 = 30 \text{ deg} \end{cases} \text{ and } \begin{cases} h_p = 300 \text{ km} \\ v_p = 10155 \text{ m/s} \end{cases} \Rightarrow \Delta V_G \approx \frac{1}{4} v_c \gamma_0^2 = 536 \text{ m/s} \Rightarrow m_f \approx 1445 \text{ kg} \quad (57)$$

The problem is then solved successively for the linear thrust model and for the bilevel thrust model with the switching date fixed to three different values, respectively 250 s, 500 s, 750 s. The Table 1 compares the solutions obtained.

| Variable | Linear model | Bilevel model | | |
|---|---|---|---|---|
| Switching date $t_1$ | – | 250 s | 500 s | 750 s |
| Thrust parameter $T_1$ | 26.467 kN | 37.942 kN | 26.339 kN | 23.378 kN |
| Thrust parameter $T_2$ | −10.976 N/s | 12.842 kN | 13.799 kN | 14.015 kN |
| Final mass $m_f$ | **1422.5 kg** | **1422.7 kg** | **1422.5 kg** | **1422.5 kg** |
| Final time $t_f$ | 1308.5 s | 1476.3 s | 1374.4 s | 1299.5 s |
| Angular range $\varphi_f$ | 69.7 deg | 79.2 deg | 73.6 deg | 69.3 deg |
| Gravity losses $\Delta V_G$ | 555 m/s | 545 m/s | 555 m/s | 558 m/s |
| AoA losses $\Delta V_T$ | 27 m/s | 37 m/s | 28 m/s | 25 m/s |

Table 1 : Optimal solutions for the linear and bilevel model

The final mass is nearly identical whatever the thrust model chosen, and it is quite close to the preliminary guess Eq. (57). The main difference between the different thrust models comes from the final time and the angular range between the initial position and the injection at the perigee. This angular range, free in the present case, has to be controlled for a transfer toward the geostationary orbit, which requires locating the apsides above the Equator. This example indicates that the perigee argument could be controlled with reduced performance loss through an adequate tuning of the thrust level.

The Figure 2 depicts the evolution of the thrust level for the linear model and for the bilevel model with the three different switching dates considered.

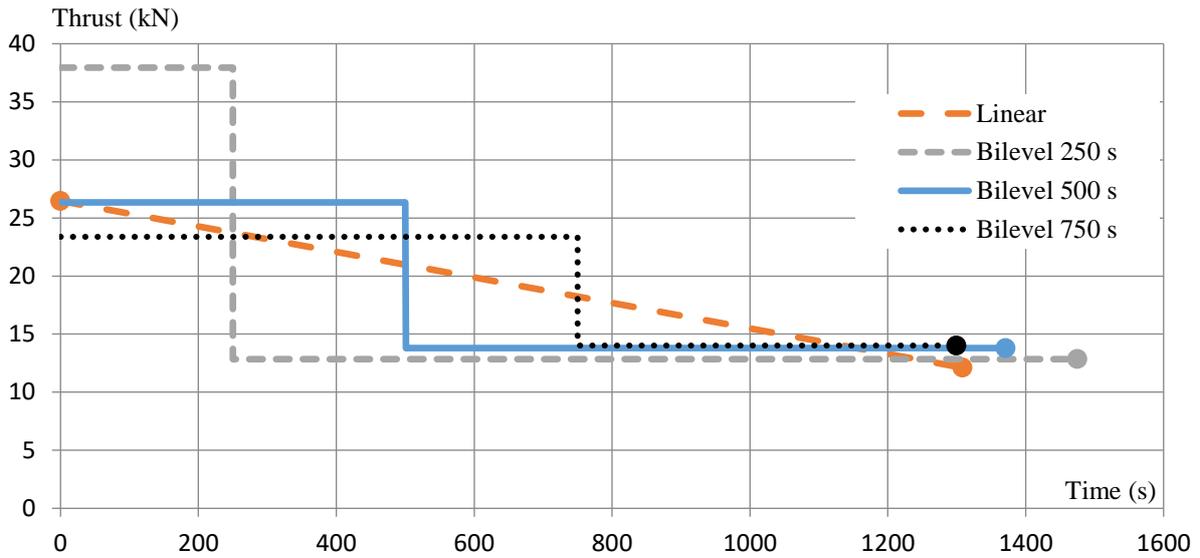

Figure 2 : Optimal thrust level

The Figure 3 and 4 plot respectively the evolution of the altitude and the angle of attack (angle between the thrust direction and the velocity) for the linear model and for the bilevel model with a switching date at 500 s.

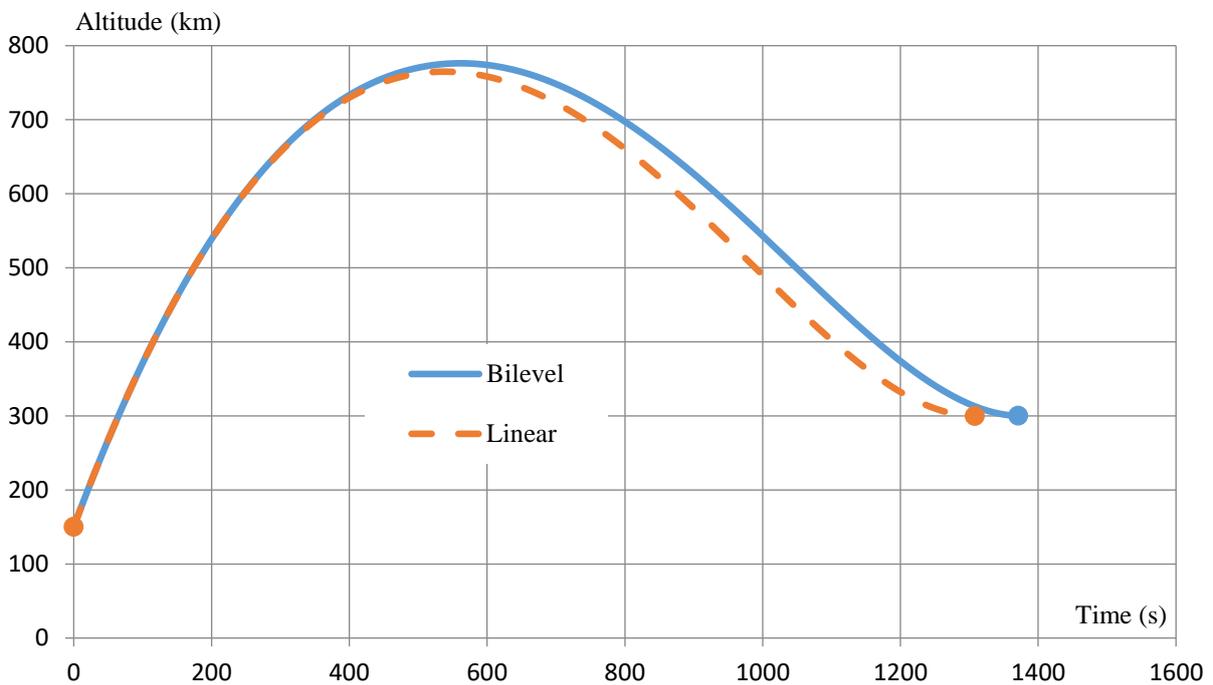

Figure 3 : Altitude

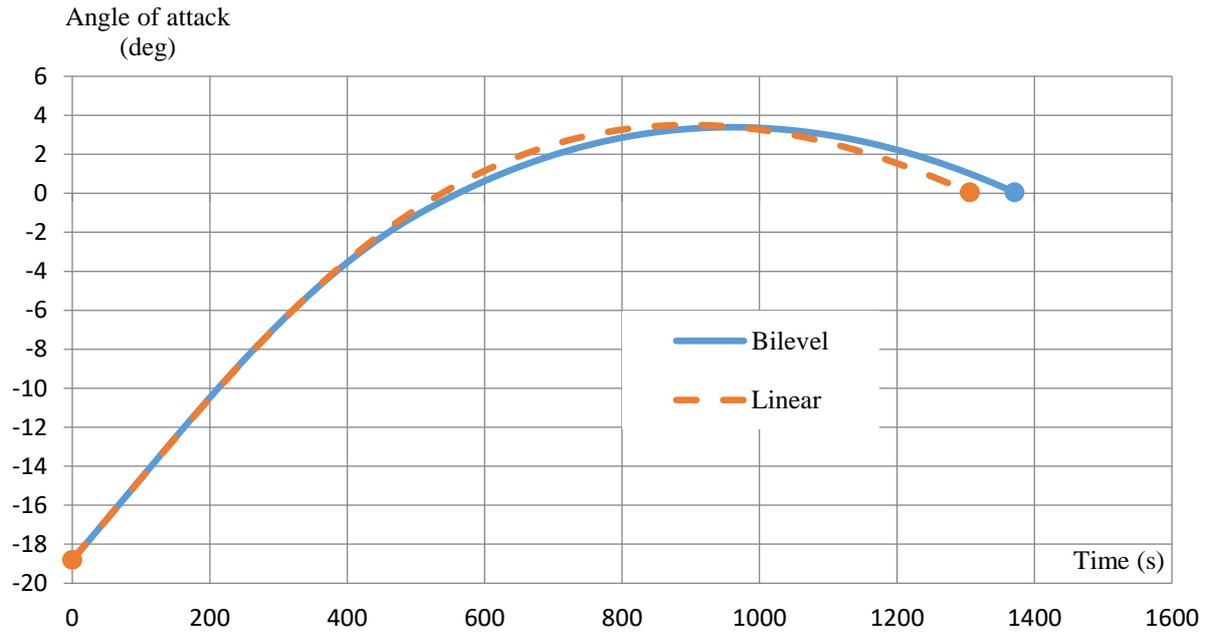

Figure 4 : Angle of attack

## 4. Conclusion

The thrust profile of a launcher upper stage is a key design parameter. The minimum-fuel single-boost planar trajectory of a launcher upper stage has been investigated in the case of a linear or a bilevel thrust profile. The profile parameters are part of the optimal control problem together with the thrust direction. Instead of complicating the optimal control problem, these additional parameters allow deriving several analytical results. First a closed-loop control is found for the thrust direction that depends only on the current kinematic conditions. The optimal control problem reduces thus to a nonlinear problem with the thrust profile parameters as only unknowns. Secondly analytical costate expressions are available that can be used to solve the transversality conditions. The propelled trajectory ends necessarily at an apside of the targeted orbit. Thirdly the direction of the trajectory final leg is downwards for a perigee injection (the most common case in practice), upwards for an apogee injection. Finally an analytical estimate of the gravity losses and of the final mass is derived.

A practical example shows that different thrust profiles lead to nearly identical values of the final mass, which was guessed with less than 2% error by the analytical estimate. The angular range depends on the thrust profile chosen, which can be useful for geostationary transfer orbits that require controlling the apside location.

These theoretical results can be used in several ways. For instance the range of reachable orbits from given initial conditions can be found by sweeping on the thrust profile parameters. Or conversely the field of extremals reaching the targeted orbit can be generated by backwards propagation from the perigee by sweeping on the thrust profile parameters.

Further work aims at extending the analysis to more complex trajectory problems. In particular some theoretical results can be hoped in the case of an initial optimized coast arc before the engine ignition, and also in the case of a three-dimensional transfer.

**Acronyms**

| | |
|---|---|
| OCP | Optimal Control Problem |
| PMP | Pontryaguin Maximum Principle |
| BVP | Boundary Value Problem |
| NLP | Nonlinear Programming |
| MDO | Multidisciplinary Optimization |
| GTO | Geostationary Transfer Orbit |